# Distributional Robust Kelly Strategy: Optimal Strategy under Uncertainty in the Long-Run

**Qingyun Sun**
School of Mathematics
Stanford University
sunqingyun1708@gmail.com

Stephen Boyd
School of Electrical Engineering
Stanford University
boyd@stanford.edu

## Abstract

In classic Kelly gambling, bets are chosen to maximize the expected log growth of wealth, under a known probability distribution. Breiman [1, 2] provides rigorous mathematical proofs that Kelly strategy maximizes the rate of asset growth (asymptotically maximal magnitude property), which is thought of as the principal justification for selecting expected logarithmic utility as the guide to portfolio selection. Despite very nice theoretical properties, the classic Kelly strategy is rarely used in practical portfolio allocation directly due to practically unavoidable uncertainty. In this paper we consider the distributional robust version of the Kelly gambling problem, in which the probability distribution is not known, but lies in a given set of possible distributions. The bet is chosen to maximize the worst-case (smallest) expected log growth among the distributions in the given set.

Computationally, this distributional robust Kelly gambling problem is convex, but in general need not be tractable. We show that it can be tractably solved in a number of useful cases when there is a finite number of outcomes with standard tools from disciplined convex programming.

Theoretically, in sequential decision making with varying distribution within a given uncertainty set, we prove that distributional robust Kelly strategy asymptotically maximizes the worst-case rate of asset growth, and dominants any other essentially different strategy by magnitude. Our results extends Breiman's theoretical result and justifies that the distributional robust Kelly strategy is the optimal strategy in the long-run for practical betting with uncertainty.

## 1 Introduction

The Classic Kelly strategy maximizes the expected logarithmic utility. It was proposed by John Kelly in a 1956 classic paper [3]. The earliest discussion of logarithmic utility dates back to 1730 in connection to Daniel Bernoulli's discussion [4] of the St. Petersburg game. In 1960 and 1961, Breiman [1, 2] proved that logarithmic utility was clearly distinguished by its optimality properties from other (essentially) different utilities as a guide to portfolio selection. The most important property is the asymptotically maximal magnitude property, which states that the Kelly strategy asymptotically maximizes the growth rate of assets and dominates the growth rate of any essentially different strategies by magnitude. As Thorp commented in [5], "This is to our mind the principal justification for selecting $\mathbf{E} \log X$ as the guide to portfolio selection."

Although the classic Kelly strategy has very nice theoretical properties, it is rarely used in practical portfolio allocation directly. The major hurdle is that in practical portfolio allocation the estimated nominal distribution of return is almost never accurate and uncertainty is unavoidable. Such uncertainty arises from multiple sources. First, the empirical nominal distribution will invariably differ from the unknown true distribution that generated the training samples, so uncertainty comes from the gap between empirical in-sample (training) distribution and out-of-sample distribution, and implementing the optimal decisions using empirical nominal distribution often leads to disappointment in out-of-sample tests. In decision analysis this phenomenon is sometimes termed the optimizer's curse. Second, uncertainty can come from a distribution shift from the model training environment and the model deployment environment, as a special but common case, in investment the major difficulty is the return time series is typically non-stationary. Third, in investment, the estimation errors from even the mean estimation and covariance estimation could be quite significant, and the further measurements of higher order distributional information beyond the first two moments of return might be uninformative due to noise in measurements.

In this paper, we identify the uncertainty problem of the classic Kelly strategy, and propose the distributional robust version of the Kelly strategy in which the probability distribution is not known, but lies in a given set of possible distributions. The distributional robust version of Kelly strategy choose bets to maximize the worst-case (smallest) expected log growth among the distributions in the given set.

Computationally, this distributional robust Kelly gambling problem is convex, but in general need not be tractable. In this work, we show that for a large class of uncertainty sets, the distributional robust Kelly problems can be transformed to tractable form that follow disciplined convex programming (DCP) rules. DCP tractable form can be solved via domain specific languages like CVXPY.

Theoretically, we extend Breiman's asymptotically maximal magnitude result and proved that distributional robust version Kelly's strategy asymptotically maximizes the worst-case rate of asset growth when the sequence of probabilities vary in the uncertainty set.

Numerically, we also tested the algorithm in a horse race gambling numerical example, and we indeed observe significant improvement of worst-case wealth growth.

**Gambling.** We consider a setting where a gambler repeatedly allocates a fraction of her wealth (assumed positive) across $n$ different bets in multiple rounds. We assume there are $n$ bets available to the gambler, who can bet any nonnegative amount on each of the bets. We let $b \in \mathbb{R}^n$ denote the bet allocation (in fraction of wealth), so $b \geq 0$ and $\mathbf{1}^T b = 1$, where $\mathbf{1}$ is the vector with all entries one. Letting $S_n$ denote the probability simplex in $\mathbb{R}^n$, we have $b \in S_n$. With bet allocation $b$, the gambler is betting $W b_i$ (in dollars) on outcome $i$, where $W > 0$ is the gambler's wealth (in dollars).

We let $r \in \mathbb{R}^n_+$ denote the random returns on the $n$ bets, with $r_i \geq 0$ the amount won by the gambler for each dollar she puts on bet $i$. With allocation $b$, the total she wins is $r^T b W$, which means her wealth increases by the (random) factor $r^T b$. We assume that the returns $r$ in different rounds are IID. We will assume that $r_n = 1$ almost surely, so $b_n$ corresponds to the fraction of wealth the gambler holds in cash; the allocation $b = e_n := (0, \dots, 0, 1)$ corresponds to not betting at all. Since her wealth is multiplied in each round by the IID random factor $r^T b$, the log of the wealth over time is therefore a random walk, with increment distribution given by the random variable $\log(r^T b)$.

**Finite outcome case.** We consider here the case where one of $K$ events occurs, *i.e.*, $r$ is supported on only $K$ points. We let $r_1, \dots, r_K$ denote the return vectors, and $\pi = (\pi_1, \dots, \pi_K) \in S_K$ the corresponding probabilities. We collect the $K$ payoff vectors into a matrix $R \in \mathbb{R}^{n \times K}$, with columns $r_1, \dots, r_K$. The vector $R^T b \in \mathbb{R}^K$ gives the wealth growth factor in the $K$ possible outcomes. The mean log growth rate is

$$G_\pi(b) = \mathbf{E}_\pi \log(r^T b) = \pi^T \log(R^T b) = \sum_{k=1}^{K} \pi_k \log(r_k^T b),$$

where the log in the middle term is applied to the vector elementwise. This is the mean drift in the log wealth random walk.

**Kelly gambling.** In a 1956 classic paper [3], John Kelly proposed to choose the allocation vector $b$ so as to maximize the mean log growth rate $G_\pi(b)$, subject to $b \geq 0$, $\mathbf{1}^T b = 1$. This method was called the Kelly criterion; since then, much work has been done on this topic [6, 5, 7, 8, 9, 10]. The mean log growth rate $G_\pi(b)$ is a concave function of $b$, so choosing $b$ is a convex optimization problem [11, 12]. It can be solved analytically in simple cases, such as when there are $K = 2$ possible outcomes. It is easily solved in other cases using standard methods and algorithms, and readily expressed in various domain specific languages (DSLs) for convex optimization such as CVX [13], CVXPY [14, 15], Convex.jl [16], or CVXR [17]. We can add additional convex constraints on $b$, which we denote as $b \in B$, with $B \subseteq S_K$ a convex set. These additional constraints preserve convexity, and therefore tractability, of the optimization problem. While Kelly did not consider additional constraints, or indeed the use of a numerical optimizer to find the optimal bet allocation vector, we still refer to the problem of maximizing $G_\pi(b)$ subject to $b \in B$ as the Kelly (gambling) problem (KP). There have been many papers exploring and extending the Kelly framework; for example, a drawdown risk constraint, that preserves convexity (hence, tractability) is described in [18]. The Bayesian version of Kelly optimal betting is described in [19]. In [20], Kelly gambling is generalized to maximize the proportion of wealth relative to the total wealth in the population.

**Distributional robust Kelly gambling.** In this paper we study a distributional robust version of Kelly gambling, in which the probability distribution $\pi$ is not known. Rather, it is known that $\pi \in \Pi$, a set of possible distributions. We define the worst-case log growth rate (under $\Pi$) as

$$G_\Pi(b) = \inf_{\pi \in \Pi} G_\pi(b).$$

This is evidently a concave function of $b$, since it is an infimum of a family of concave functions of $b$, *i.e.*, $G_\pi(b)$ for $\pi \in \Pi$. The *distributional robust Kelly problem* (DRKP) is to choose $b \in B$ to maximize $G_\Pi(b)$,

$$\begin{array}{ll} \text{maximize} & \inf_{\pi \in \Pi} \mathbf{E}_\pi \log(r^T b) \\ \text{subject to} & b \in B. \end{array}$$

This is in principle a convex optimization problem, specifically a distributional robust problem; but such problems in general need not be tractable, as discussed in [21, 22, 23] The purpose of this paper is to show how the DRKP can be tractably solved for some useful probability sets $\Pi$ via disciplined convex programming. In this paper we call an optimization problem "DCP tractable" in a strict and specific sense when the optimization problem is a disciplined convex problem, so that it can be solved via domain specific languages like CVXPY.

**Related work on uncertainty aversion.** In decision theory and economics, there are two important concepts, risk and uncertainty. Risk is about the situation when a probability can be assigned to each possible outcome of a situation. Uncertainty is about the situation when the probabilities of outcomes are unknown. Uncertainty aversion, also called ambiguity aversion, is a preference for known risks over unknown risks. Uncertainty aversion provides a behavioral foundation for maximizing the utility under the worst of a set of probability measures; see [24, 25, 26, 27] for more detailed discussion. The Kelly problem addresses risk; the distributional robust Kelly problem is a natural extension that considers uncertainty aversion.

**Related work on distributional robust optimization.** Distributional robust optimization is a well studied topic. Previous work on distribution robust optimization studied finite-dimensional parametrization for probability sets including moments, support or directional deviations constraints in [28, 29, 30, 31, 32, 33]. Beyond finite-dimensional parametrization of the probability set, researchers have also studied non-parametric distances for probability measure, like $f$-divergences (*e.g.*, Kullback-Leibler divergences) [34, 35, 36, 37, 38] and Wasserstein distances [39, 40, 41, 42].

**Contribution**

- We propose distributional robust Kelly strategy to account for uncertainty in practical investment, where we maximize the worst-case (smallest) expected log growth among the distributions in the given set. Theoretically, we proved that distributional robust version Kelly strategy asymptotically maximizes the worst-case rate of asset growth when the sequence of probabilities vary in the uncertainty set and dominant any other essentially different strategy by magnitude.
- For a large class of uncertainty sets including polyhedral sets, ellipsodial sets, $f$-divergence ball, Wasserstein ball and uncertainty set with estimated mean and covariance, we concretely derived the disciplined convex programming (DCP) forms of distributional robust Kelly problems and provided concrete software implementations using CVXPY, with a simple horse gamble example to numerically verify that distributional robust Kelly strategies indeed lead to a better worst-case wealth growth, and also lead to more diverse bet vector in this example.

## 2 DCP tractable forms for distributional robust Kelly strategy

In this section we show how to formulate DRKP as a DCP tractable convex optimization problem for a variety of distribution sets. The key is to derive a DCP tractable description of the worst-case log growth $G_\Pi(b)$. We use duality to express $G_\Pi(b)$ as the value of a convex maximization problem, which allows us to solve DRKP as one convex problem.

**Polyhedron defined by linear inequalities and equalities.**

**Theorem 1.** *For polyhedron uncertainty set given by a finite set of linear inequalities and equalities,*

$$\Pi = \{\pi \in S_K \mid A_0\pi = d_0,\ A_1\pi \leq d_1\},$$

*where $A_0 \in \mathbb{R}^{m_0\times K}$, $b_0 \in \mathbb{R}^{m_0}$, $A_1 \in \mathbb{R}^{m_1\times K}$, $b_1 \in \mathbb{R}^{m_1}$, the distributional robust Kelly problem is*

$$\begin{array}{ll} \textit{maximize} & \min(\log(R^T b) + A_0^T\mu + A_1^T\lambda) - d_0^T\mu - d_1^T\lambda \\ \textit{subject to} & b \in B, \quad \lambda \geq 0, \end{array}$$

*with variables $b, \mu, \lambda$. The problem follows the disciplined convex programming (DCP) rules.*

**Box uncertainty set.** As a commonly used special case of polyhedron, we consider box.

**Theorem 2.** *for box uncertainty set*

$$\Pi = \{\pi \in S_K \mid |\pi - \pi^{\rm nom}| \leq \rho\},$$

*, where $\pi^{\rm nom} \in S_K$ is the nominal distribution, and $\rho \in \mathbb{R}^n_+$ is a vector of radii, (The inequality $|\pi - \pi^{\rm nom}| \leq \rho$ is interpreted elementwise), the distributional robust Kelly problem is*

$$\begin{array}{ll} \textit{maximize} & \min(\log(R^T b) + \lambda) - (\pi^{\rm nom})^T\lambda - \rho^T|\lambda| \\ \textit{subject to} & b \in B, \end{array}$$

*with variables $b, \lambda$. The problem follows the disciplined convex programming (DCP) rules.*

**Ellipsoidal uncertainty set** Here we consider the case when $\Pi$ is the inverse image of a $p$-norm ball, with $p \geq 1$, under an affine mapping. As usual we define $q$ by $1/p + 1/q = 1$. This includes an ellipsoid (and indeed the box set described above) as a special case.

**Theorem 3.** *For ellipsoidal uncertainty set*

$$\Pi = \{\pi \in S_K \mid \|W^{-1}(\pi - \pi^{\rm nom})\|_p \leq 1\},$$

*where $W$ is a nonsingular matrix, the distributional robust Kelly problem is*

$$\begin{array}{ll} \textit{maximize} & (\pi^{\rm nom})^T(u) - \|W^T(u - \mu 1)\|_q \\ \textit{subject to} & u \leq \log(R^T b), \\ & b \in B, \end{array}$$

*with variables $b, u, \mu$. The problem follows the disciplined convex programming (DCP) rules.*

This DRKP problem follows DCP rule because $\pi^{{\rm nom},T}u$ is a linear function of $u$, $-\|W^T(u-\mu 1)\|_q$ is a concave function of $u$ and $\mu$ for $q \geq 1$, and $\log(R^T b)$ is a concave function of $b$, hence $u \leq \log(R^T b)$ is a concave constraint. Therefore, this DRKP problem is maximizing a concave objective concave constraint problem that follows DCP rule.

The proof of the theorem is based on the Lagrangian duality and Hölder equality for the $p$-norm:

$$\sup_{\|z\|_p\leq 1} z^T W^T x = \|W^T x\|_q,$$

**Divergence based distribution set** Let $\pi_1, \pi_2 \in S_K$ be two distributions. For a convex function $f : \mathbb{R}_+ \to \mathbb{R}$ with $f(1) = 0$, the $f$-divergence of $\pi_1$ from $\pi_2$ is defined as

$$D_f(\pi_1\|\pi_2) = \pi_2^T f(\pi_1/\pi_2),$$

where the ratio is meant elementwise. Recall that the Fenchel conjugate of $f$ is $f^*(s) = \sup_{t\geq 0}(ts - f(t))$.

**Theorem 4.** *For $f$-divergence ball uncertainty set*

$$\Pi = \{\pi \in S_K \mid D_f(\pi \| \pi^{\text{nom}}) \leq \epsilon\},$$

*where $\epsilon > 0$ is a given value, the distributional robust Kelly problem is*

$$\begin{array}{ll} \textit{maximize} & -(\pi^{\text{nom}})^T w - \epsilon\lambda - \gamma \\ \textit{subject to} & w \geq \lambda f^*(\frac{z}{\lambda}) \\ & z \geq -\log(R^T b) - \gamma \\ & \lambda \geq 0, \quad b \in B, \end{array}$$

*with variables $b, \gamma, \lambda, w, z$.*

*The problem follows the disciplined convex programming (DCP) rules.*

Here

$$\lambda f^*(\frac{z}{\lambda}) = (\lambda f)^*(z) = \sup_{t \geq 0}(tz - \lambda f(t)),$$

is the perspective function of the non-decreasing convex function $f^*(z)$, so it is also a convex function that is non-decreasing in $z$. Additionally, $-\log(R^T b) - \gamma$ is a convex function of $b$ and $\gamma$; then from the DCP composition rule, we know this form of DRKP is convex. Concrete examples of $f-$divergence function and their Fenchel conjugate funuctions are provided in supplementary material for the convenience of readers.

**Wasserstein distance uncertainty set with finite support** When $\pi$ and $\pi^{\text{nom}}$ both have finite supports, the Wasserstein distance $D_c(\pi, \pi^{\text{nom}})$ with cost $c \in \mathbb{R}_+^{K \times K^{\text{nom}}}$ is defined as the optimal value of the problem

$$\begin{array}{ll} \text{minimize} & \sum_{i,j} Q_{ij} c_{ij} \\ \text{subject to} & Q\mathbf{1} = \pi, \quad Q^T\mathbf{1} = \pi^{\text{nom}}, \quad Q \geq 0, \end{array}$$

with variable $Q$. The Wasserstein distance has several other names, including Monge-Kantorovich, earth-mover, or optimal transport distance [39, 40, 41, 42]. The Wasserstein uncertainty set would be better to prepare for black swans events comparing to $f$-divergence uncertainty set, as it does not require $\pi$ to be absolutely continuous with respect to the nominal distribution.

**Theorem 5.** *For Wasserstein distance ball uncertainty set*

$$\Pi = \{\pi \in S_K \mid D_c(\pi, \pi^{\text{nom}}) \leq s\},$$

*with $s > 0$. the distributional robust Kelly problem is*

$$\begin{array}{ll} \textit{maximize} & \left(\sum_j \pi_j^{\text{nom}} \min_i(\log(R^T b)_i + \lambda c_{ij}) - s\lambda\right) \\ \textit{subject to} & b \in B, \quad \lambda \geq 0, \end{array}$$

*where $\lambda \in \mathbb{R}_+$ is the dual variable.*

*The problem follows the disciplined convex programming (DCP) rules.*

The problem follows the disciplined convex programming (DCP) rules, because $\log(R^T b)_i + \lambda c_{ij}$ is a concave function of $b$ and $\lambda$, therefore $\min_i(\log(R^T b)_i + \lambda c_{ij})$ is a concave function of $b$ and $\lambda$, then the entire objective is a concave function of $b$ and $\lambda$; and the constraint $b \in B$ and $\lambda \geq 0$ also follows DCP rule. We comment that although we allow $\pi$ to have different support from $\pi^{\text{nom}}$, we still consider the simple setting of finite event space for $\pi$ for technical clarity. The extension to the general setting for different norm form of cost $c$ could be find in [43]. Computing the Wasserstein distance between two discrete distributions can be converted to solving a tractable linear program that is susceptible to the network simplex algorithm, dual ascent methods, or specialized auction algorithms [44, 45, 46]. Efficient approximation schemes can be found in the survey [47] of algorithms for the finite-dimensional transportation problem. However, as soon as at least one of the two involved distributions is not discrete, the Wasserstein distance can no longer be evaluated in polynomial time.

**Uncertainty set with estimated mean and covariance matrix of return** For the application in stock investment, typically quantitative investors are only able to obtain estimated mean and covariance matrix, with error bounds on the estimation errors. Following the original paper by Delage and Ye [28] and the review paper [48], we consider the following uncertainty set.

**Theorem 6.** *For uncertainty set with estimator $\mu_0$ and $\Sigma_0$ for mean and covariance matrix of random vector $r \in I\!R^n$,*

$$\Pi = \{\pi \in S_K \mid (\mathbf{E}_\pi r - \mu_0)^T \Sigma_0^{-1} (\mathbf{E}_\pi r - \mu_0) \le \varrho_1, \quad \mathbf{E}_\pi[(r-\mu_0)(r-\mu_0)^T] \preccurlyeq \varrho_2\},$$

*if there exists $\pi_0 \in \Pi$ such that $\mathbf{E}_{\pi_0} r = \mu_0$, $\mathbf{E}_{\pi_0}(r-\mu_0)(r-\mu_0)^T = \Sigma_0$, then the distributional robust Kelly problem has the same optimal value as the following SDP problem,*

$$\begin{array}{ll} \textit{minimize} & u_1 + u_2 \\ \textit{subject to} & u_1 \ge -\log(r_i^T b) - r_i^\top \boldsymbol{Y} r_i - r_i^\top \boldsymbol{y}, \ \forall i = 1,\dots,K, \\ & u_2 \ge (\varrho_2 \boldsymbol{\Sigma}_0 + \boldsymbol{\mu}_0 \boldsymbol{\mu}_0^\top) \bullet \boldsymbol{Y} + \boldsymbol{\mu}_0^\top \boldsymbol{y} + \sqrt{\varrho_1} \|\boldsymbol{\Sigma}_0^{\frac{1}{2}}(\boldsymbol{y} + 2\boldsymbol{Y}\boldsymbol{\mu}_0)\|, \\ & \boldsymbol{Y} \succcurlyeq 0, \\ & b \in B, \end{array}$$

*where $u_1, u_2 \in I\!R$, $\boldsymbol{Y} \in I\!R^{n\times n}$ and $\boldsymbol{y} \in I\!R^n$ are auxiliary variables.*

*The problem is a SDP problem, therefore follows the disciplined convex programming (DCP) rules.*

## 3 Theoretical properties of distributional robust Kelly strategy

In the following, we will extend Brieman's classical optimality property of the Kelly strategy to distributional robust Kelly strategy. We show that for sequential decision making problems under uncertainty, distributional robust Kelly strategy is

Consider a sequential gambling setting with a fixed uncertainty set $\Pi$. For the $N$th period random return vector $r_N$ has distribution $\pi_N$. $\pi_N$ could vary freely for different $N$ and the only condition is that $\pi_N \in \Pi$ for any $N$. In this theoretical section, we consider the more general setting where the event (outcome) space is not necessarily finite, we only assume that return is bounded $r_{N,i} \in [0, R_M]$ for all $i, N$. let $b_N$ be the proportional betting vector for the $N$-th period such that $b_N \ge 0$, $1^T b_N = 1$, denote the (sequential) strategy $\Lambda = (b_1, \dots, b_N, \dots)$. This strategy's wealth growth at the $N$-th period is $V_N := r_N^T b_N$. This strategy's accumulated wealth growth is $S_N$ , defined recursively as $S_0 = 1$, $S_{N+1} = S_N V_N$. Let $O_{N-1}$ be the outcomes during the first $N-1$ investment periods. Let $V_N^\star := r_N^T b_N^\star$. $S_N^\star = S_{N-1}^\star V_{N-1}^\star$ where $b_N^\star$ is the distributional robust bet at the $N$-th period that maximizes

$$\inf_{\pi_N \in \Pi} \mathbf{E}_{\pi_N}[\log(r^T b) \mid O_{N-1}]$$

**Theorem 7.** *For any strategy $\Lambda$ leading to the fortune $S_N$, $\lim_N \frac{S_N}{S_N^\star}$ exists almost surely and $\inf_{\pi_1,\dots,\pi_N \in \Pi} \mathbf{E} \lim_N \frac{S_N}{S_N^\star} \le 1$.*

The critical proof idea is to leverage the maximizing property of the distributional robust bet $b_N^\star$, using superlinear property of the operator $\mathcal{E} := \inf_{\pi_N \in \Pi} \mathbf{E}_{\pi_N}$: $\inf_{\pi_N \in \Pi} \mathbf{E}_{\pi_N}(S^1 + S^2) \ge \inf_{\pi_N \in \Pi} \mathbf{E}_{\pi_N}(S^1) + \inf_{\pi_N \in \Pi} \mathbf{E}_{\pi_N}(S^2)$, eventually using Fatou's lemma to switch the order of limit and expectation.

*Proof sketch .* Here we highlight some of the key technical points, with the full proof in the supplementary material. We have

$$\inf_{\pi_N \in \Pi} \mathbf{E}_{\pi_N}[\frac{S_N}{S_N^\star} \mid O_{N-1}] = \inf_{\pi_N \in \Pi} \mathbf{E}_{\pi_N}[\frac{V_N}{V_N^\star} \mid O_{N-1}] \frac{S_{N-1}}{S_{N-1}^\star}$$

with

$$\inf_{\pi_1,\dots,\pi_{N-1},\pi_N \in \Pi} \mathbf{E}_{\pi_1,\dots,\pi_{N-1},\pi_N} \lim_N \frac{S_N}{S_N^\star} \le \inf_{\pi_1,\dots,\pi_{N-1},\pi_{N-1} \in \Pi} \mathbf{E}_{\pi_1,\dots,\pi_{N-1}} \lim_N \frac{S_{N-1}}{S_{N-1}^\star} \le \frac{S_0}{S_0^\star} = 1.$$

To prove the theorem, we only need prove that for any $b_N$,

$$\inf_{\pi_N \in \Pi} \mathbf{E}_{\pi_N}(\frac{V_N}{V_N^\star} \mid O_{N-1}) \leq 1$$

Now, for any $\epsilon > 0$, by the maximizing property of the distributional robust bet $b_N^\star$, we have

$$\inf_{\pi_N \in \Pi} \mathbf{E}(\epsilon \log(V_N) + (1-\epsilon)\log(V_N^\star) \mid O_{N-1}) \leq \inf_{\pi_N \in \Pi} \mathbf{E}(\log(V_N^\star) \mid O_{N-1})$$

Rewrite the left side, using superlinear property of $\mathcal{E} := \inf_{\pi_N \in \Pi} \mathbf{E}_{\pi_N}$, we have

$$\inf_{\pi_N \in \Pi} \mathbf{E}[\frac{1}{\epsilon}\log(1 + \frac{\epsilon}{1-\epsilon}\frac{V_N}{V_N^\star}) \mid O_{N-1}] \leq \frac{1}{\epsilon}\log(\frac{1}{1-\epsilon})$$

taking lower limit $\epsilon \to 0^+$, since $\frac{V_N}{V_N^\star}$) are bounded, from Fatou's lemma, we have

$$\inf_{\pi_N \in \Pi} \mathbf{E}[\frac{V_N}{V_N^\star}) \mid O_{N-1}] = \inf_{\pi_N \in \Pi} \mathbf{E}[\underline{\lim}_{\epsilon \to 0^+}\frac{1}{\epsilon}\log(1 + \frac{\epsilon}{1-\epsilon}\frac{V_N}{V_N^\star}) \mid O_{N-1}] \leq \underline{\lim}_{\epsilon \to 0^+}\frac{1}{\epsilon}\log(\frac{1}{1-\epsilon}) = 1$$

□

We call $\Lambda$ a nonterminating strategy if there are no values of $b_s$ such that $V_s = r_s^T b_s = 0$ for any $s$.

**Theorem 8.** *If $\Lambda$ is a nonterminating strategy, the set $\overline{\lim}_{N\to\infty}\frac{S_N}{S_N^\star} = 0$ is almost surely equal to the set on which $\sum_{N=1}^{\infty}[\inf_{\pi_N \in \Pi} \mathbf{E}_{\pi_N}[\log(V_N^\star) - \log(V_N) \mid O_{N-1}]] = \infty$.*

The critical proof idea is to combine previous theorem and a generalized martingale convergence theorem on the sequence

$$\frac{S_N^\star}{S_N} - \inf_{\pi_N \in \Pi} \mathbf{E}_{\pi_s}[\frac{S_N^\star}{S_N} \mid O_{N-1}] = \sum_{s=1}^{N}\{\log(V_s^\star) - \log(V_s) - \inf_{\pi_s \in \Pi} \mathbf{E}[\log(V_s^\star) - \log(V_s) \mid O_{s-1}]\}$$

using non-linear expectation theory developed by [49, 50] on the non-linear expectation operator $\mathcal{E} := \inf_{\pi_N \in \Pi} \mathbf{E}_{\pi_N}$.

We comment that if the two strategies $\Lambda$ and $\Lambda^\star$ satisfy the condition $\sum_{N=1}^{\infty}[\inf_{\pi_N \in \Pi} \mathbf{E}_{\pi_N}[\log(V_N^\star) - \log(V_N) \mid O_{N-1}]] = \infty$, then we call $\Lambda$ and $\Lambda^\star$ "essentially different" strategies under uncertainty. The two theorems proved in this section could be stated as:

*In sequential decision making problem under uncertainty, when the sequence of distributions of return vary in a given uncertainty set, distributional robust Kelly strategy asymptotically maximizes the worst-case rate of asset growth, and dominants any other essentially different strategy by magnitude.*

# 4 Numerical example

In this section we illustrate distributional robust Kelly gambling with a simple horse racing example. Our example is a simple horse race with $n$ horses, with bets placed on each horse placing, *i.e.*, coming in first or second. There are thus $K = n(n-1)/2$ outcomes (indexed as $j, k$ with $j < k \leq n$), and $n$ bets (one for each horse to place). We consider two simple uncertainty sets, the box set with parameter $\eta$ and $\ell_2$ ball set with radius $c$. We first describe the nominal distribution of outcomes $\pi^{\text{nom}}$. We model the speed of the horses as independent random variables, with the fastest and second fastest horses placing. With this model, $\pi^{\text{nom}}$ is entirely described by the probability that horse $i$ comes in first, we which denote $\beta_i$. For $j < k$, we have $\pi_{jk}^{\text{nom}} = P(\text{horse } j \text{ and } k \text{ are the first two}) = \beta_j\beta_k(\frac{1}{1-\beta_i} + \frac{1}{1-\beta_j})$.

For the return matrix, we use parimutuel betting, with the fraction of bets on each horse equal to $\beta_i$, the probability that it will win (under the nominal probability distribution). The return matrix $R \in \mathbb{R}^{n\times K}$ then has the form (we index the columns (outcomes) by the pair $jk$, with $j < k$)

$$R_{i,jk} = \begin{cases} \frac{n}{1+\beta_j/\beta_k} & i = j \\ \frac{n}{1+\beta_k/\beta_j} & i = k \\ 0 & i \notin \{j,k\}, \end{cases}$$

| **Growth rate** | $b^K$ | $b^{RK}$: box with $\eta = 0.26$ | $b^{RK}$: ball with $c = 0.016$ |
|---|---|---|---|
| $\pi^{\text{nominal}}$ | 4.3% | 2.2% | 2.2% |
| $\pi^{\text{worst}}$ | −2.2% | 0.7% | 0.4% |

**Table 1:** For box uncertainty set with $\eta = 0.26$ and ball uncertainty set with $c = 0.016$, we compare the growth rate and worst-case growth rate for the Kelly optimal and the distributional robust Kelly optimal bets.

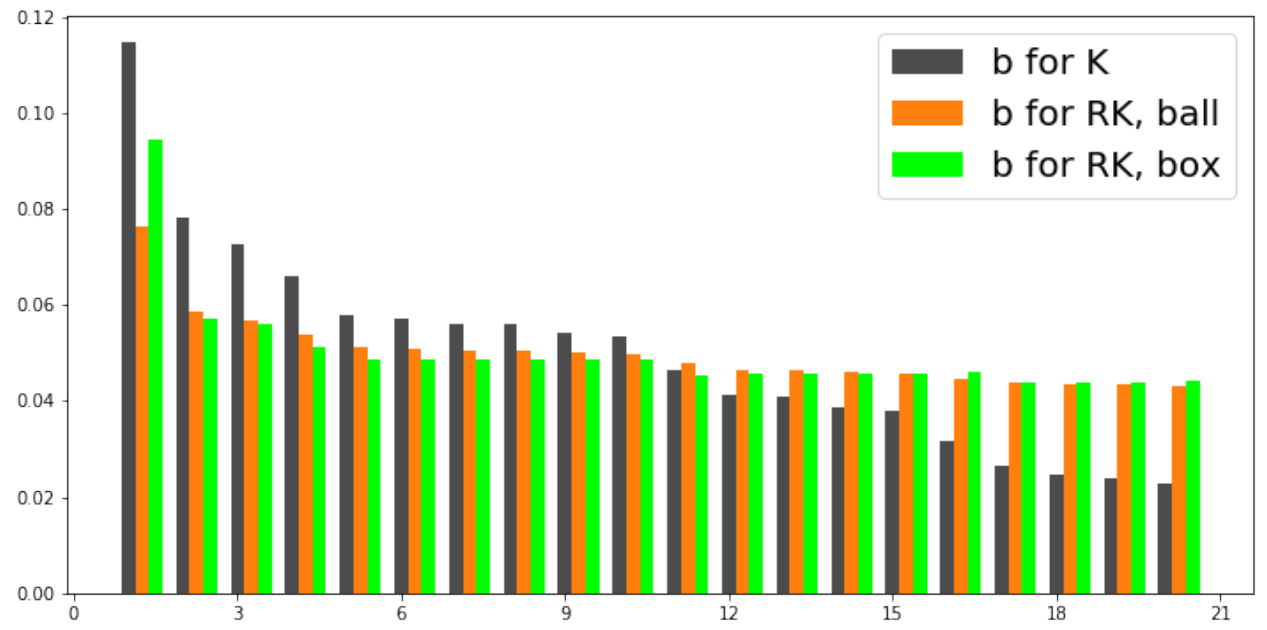

**Figure 1:** The Kelly optimal bets $b^K$ for the nominal distribution, and the distributional robust optimal bets for the box and ball uncertainty sets, ordered by the descending order of $b^K$.

First, we show growth rate and worst-case growth rate for the Kelly optimal and the distributional robust Kelly optimal bets under two uncertainty sets. In table 1, we show the comparison for box uncertainty set with $\eta = 0.26$ and for ball uncertainty set with $c = 0.016$. The two parameters are chosen so that the worst case growth of Kelly bets for both uncertainty sets are $-2.2\%$. In particular, using standard Kelly betting, we lose money (when the distribution is chosen as the worst one for the Kelly bets). We can see that, as expected, the Kelly optimal bet has higher log growth under the nominal distribution, and the distributional robust Kelly bet has better worst-case log growth. We see that the worst-case growth of the distributional robust Kelly bet is significantly better than the worst-case growth of the nominal Kelly optimal bet. In particular, with robust Kelly betting, we make money, even when the worst distribution is chosen. The nominal Kelly optimal bet $b^{\text{K}}$ and the distributional robust Kelly bet $b^{\text{RK}}$ for both uncertainty sets in figure 11. For each of our bets $b^{\text{K}}$ and $b^{\text{RK}}$ shown above, we find a corresponding worst case distribution, denoted $\pi^{\text{wc,K}}$ and $\pi^{\text{wc,RK}}$, which minimize $G_\pi(b)$ over $\pi \in \Pi$. These distributions, shown for box uncertainty set and ball uncertainty set in figure 3, achieve the corresponding worst-case log growth for the two bet vectors. Finally, in figure 2, we compare the expected wealth logarithmic growth rate as we increase the size of the uncertainty sets. For the box uncertainty set we choose $\eta \in [0, 0.3]$, and for the ball uncertainty set we choose $c \in [0, 0.02]$, we look at the expected growth for both the kelly bet $b^K$ and the distributional robust kelly bet $b^{RK}$ under both the nominal probability $\pi^{\text{nom}}$ and the worst case probability $\pi^{\text{worst}}$.

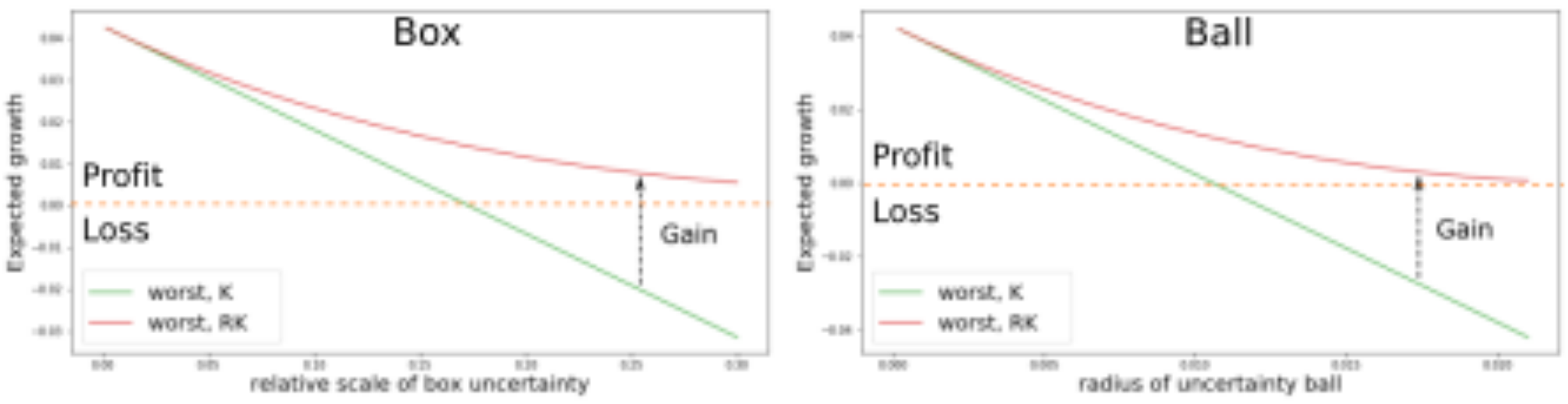

**Figure 2:** Plots of the expected growth under nominal distribution and worst-case distribution under the box and ball uncertainty set family. The blue, green, orange, red line are $\pi^{\text{nom},T} \log(R^T b^K)$, $\pi_\eta^{\text{worst},T} \log(R^T b^K)$, $\pi^{\text{nom},T} \log(R^T b_\eta^{RK})$, $\pi_\eta^{\text{wc},T} \log(R^T b_\eta^{RK})$. For the box uncertainty set we choose $\eta \in [0, 0.3]$, and for the ball uncertainty set we choose $c \in [0, 0.02]$.

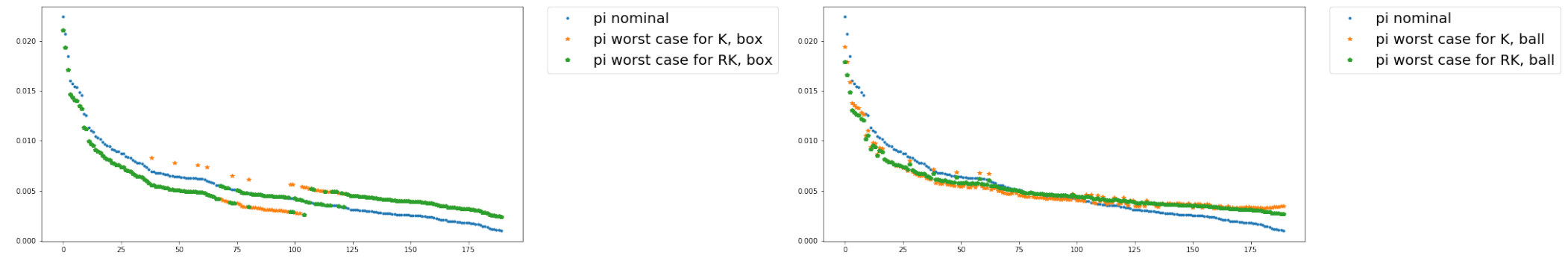

**Figure 3:** For box uncertainty set with $\eta = 0.26$ and ball uncertainty set with $c = 0.016$, we show the nominal distribution $\pi^{\text{nom}}$ (sorted) and the two worst-case distributions $\pi^{\text{wc,K}}$ and $\pi^{\text{wc,RK}}$.

## 5 Discussion about how to choose uncertainty set

**Uncertainty set estimation via chance constraints** To use distributional robust Kelly strategy in practical, one of the limitation is that it could be hard to estimate a high quality uncertainty set and the choice of uncertainty set depends on domain knowledge. There are two major considerations in the choice of the uncertainty set, first is tractability, which is discussed through the lens of DCP tractable form in this paper, second is efficiency, or the trade-off between coverage and tightness, i.e. whether it is too conservative and how well it reflects the actual variability in our problem. We make some comment on the selection advice for uncertainty set. From the relation between chance constraint and uncertainty set, as discussed in section 6 of [51, 52], we could use probabilistic tools to efficiently represent uncertainty sets through chance constraint such as drawdown control, including value at Risk, conditional value at risk, tail bound from moments[18]. In quantitative investment practice, due to the noisy nature of stocks or futures market data, the scarcity of effective observations, the return non-stationarity, the uncertainty set presented in theorem 6 would be a good starting point. Besides the long-run optimality of growth rate proved here, an interesting question is to study the finite time property of the distributional robust Kelly strategy and compare both the achieved growth rate in finite time and the probability of having a given level of drawdown with the distributional robust version of mean-variance strategy.

**Uncertainty set estimation via conformal prediction** Recent uncertainty quantification methodology development on conformal prediction [53, 54, 55] has provided statistical tools to generate set-valued predictions for black-box predictors with rigorous error control and formal finite-sample coverage guarantee, as shown in [56, 57]. As a future direction, we will explore the usage of conformal prediction to construct computational tractable uncertainty set from data for investment.

**Uncertainty set tuning via bi-level optimization and differentiation through solution of convex problem** To automate the tuning of the parameters for the uncertainty set, we could use bi-level optimization [58] to learn uncertainty set by back-propagation over the parameter $\theta$. We could set up a higher-level objective $L(\theta)$ to tune the parameters $\theta$. For example, given out of samples observations that are not accessible when we fit the mean the functional form of $\Pi$, we could define an out-of-sample Kelly loss as the higher-level objective:

$$b^\star(\theta) = \arg\max_{b \in B} \quad \min_{\pi \in \Pi(\theta)} \mathbf{E}_\pi \log(r^T b).$$

$$\max_\theta \quad L(\theta) = \frac{1}{N} \sum_{i=1}^{N} \log(r_i^T b^\star(\theta))$$

Our method relies on recently developed methods[58] that can efficiently evaluate the derivative of the solution of a disciplined convex optimization problem with respect to its parameters. Using previous result in this paper, this distributional robust optimization problem can be transformed into a disciplined convex optimization, which allows automated differentiation with respect to the solution map $b^\star(\theta)$. We have provided CVXPY code for the uncertainty sets and their corresponding tuning/learning example in PyTorch in the supplementary material, which would allow users to learn uncertainty sets through the recently developed software framework to embed our problem into differentiable programming framework to learn the uncertainty set.

# 6 Supplementary Material: Concrete examples of divergences functions and their Fenchel conjugates

We remark that there is a one-parameter family of $f$-divergences generated by the $\alpha$-function with $\alpha \in I\!R$, where we can define the generalization of natural logarithm by

$$\log_\alpha(t) = \frac{t^{\alpha-1} - 1}{\alpha - 1}.$$

For $\alpha \neq 1$, it is a power function, for $\alpha \to 1$, it is converging to the natural logarithm. Now if we assume $f_\alpha(1) = 0$ and $f'_\alpha(t) = \log_\alpha(t)$, then we have

$$f_\alpha(t) = \frac{t^\alpha - 1 - \alpha(t-1)}{\alpha(\alpha - 1)}.$$

The Fenchel conjugate is

$$f^*_\alpha(s) = \frac{1}{\alpha}((1 + (\alpha - 1)s)^{\frac{\alpha}{\alpha-1}} - 1).$$

We now show some more specific examples of $f$-divergences; for a more detailed discussion see [37].

- *KL-divergence.* With $f(t) = t\log(t) - t + 1$, we obtain the KL-divergence. We have $f^*(s) = \exp(s) - 1$. This corresponds to $\alpha = 1$.
- *Reverse KL-divergence.* With $f(t) = -\log(t) + t - 1$, the $f$-divergence is the reverse KL-divergence. We have $f^*(s) = -\log(1 - s)$ for $s < 1$. This corresponds to $\alpha = 0$.
- *Pearson $\chi^2$-divergence.* With $f(t) = \frac{1}{2}(t-1)^2$, we obtain the Pearson $\chi^2$-divergence. We have $f^*(s) = \frac{1}{2}(s+1)^2 - \frac{1}{2}$, $s > -11$. This corresponds to $\alpha = 2$.
- *Neyman $\chi^2$-divergence.* With $f(t) = \frac{1}{2t}(t-1)^2$, we obtain the Neyman $\chi^2$-divergence. We have $f^*(s) = 1 - \sqrt{1 - 2s}$, $s < \frac{1}{2}$. This corresponds to $\alpha = -1$.
- *Hellinger-divergence.* With $f(t) = 2(\sqrt{t} - 1)^2$, we obtain the Hellinger-divergence. We have $f^*(s) = \frac{2s}{2-s}$, $s < 2$. This corresponds to $\alpha = -1$.
- *Total variation distance.* With $f(t) = |t - 1|$, the $f$-divergence is the total variation distance. We have $f^*(s) = -1$ for $s \leq -1$ and $f^*(s) = s$ for $-1 \leq s \leq 1$.

# 7 Supplementary Material: Proof of theorems in section 2

**Proof of Theorem** 1

*Proof of theorem* 1. The worst-case log growth rate $G_\Pi(b)$ is given by the optimal value of the linear program (LP)

$$\begin{array}{ll} \text{minimize} & \pi^T \log(R^T b) \\ \text{subject to} & \mathbf{1}^T \pi = 1, \quad \pi \geq 0, \\ & A_0 \pi = d_0, \quad A_1 \pi \leq d_1, \end{array} \tag{1}$$

with variable $\pi$.

We form a dual of this problem, working with the constraints $A_0\pi = d_0$, $A_1\pi \leq d_1$; we keep the simplex constraints $\pi \geq 0$, $\mathbf{1}^T\pi = 1$ as an indicator function $I_S(\pi)$ in the objective. The Lagrangian is

$$\begin{aligned} L(\nu, \lambda, \pi) = \quad & \pi^T \log(R^T b) + \nu^T(A_0\pi - d_0) + \\ & \lambda^T(A_1\pi - d_1) + I_S(\pi), \end{aligned}$$

where $\nu \in I\!R^{m_0}$ and $\lambda \in I\!R^{m_1}$ are the dual variables, with $\lambda \geq 0$. Minimizing over $\pi$ we obtain the dual function,

$$\begin{aligned} g(\nu, \lambda) \quad & = \inf_{\pi \in S_K} L(\nu, \lambda, \pi) \\ & = \min(\log(R^T b) + A_0^T \nu + A_1^T \lambda) - d_0^T \mu - d_1^T \lambda, \end{aligned}$$

where the min of a vector is the minimum of its entries. The dual problem associated with (1) is then

$$\begin{array}{ll} \text{maximize} & \min(\log(R^T b) + A_0^T \mu + A_1^T \lambda) - d_0^T \mu - d_1^T \lambda \\ \text{subject to} & \lambda \geq 0, \end{array}$$

with variables $\mu, \lambda$. Using Slater's condition for simplex, strong duality can easily be verified. Therefore, this dual problem has the same optimal value as (1), *i.e.*,

$$\begin{array}{rl} G_\Pi(b) = \sup_{\nu,\lambda\geq 0}( & \min(\log(R^T b) + A_0^T \mu + A_1^T \lambda) \\ & -d_0^T \mu - d_1^T \lambda). \end{array}$$

Using this expression for $G_\Pi(b)$, the DRKP becomes

$$\begin{array}{ll} \text{maximize} & \min(\log(R^T b) + A_0^T \mu + A_1^T \lambda) - d_0^T \mu - d_1^T \lambda \\ \text{subject to} & b \in B, \quad \lambda \geq 0, \end{array}$$

with variables $b, \mu, \lambda$. □

**Proof of Theorem 2**

*Proof of theorem* 2. Using the general method above, expressing the limits as $A_1 \pi \leq d_1$ with

$$A_1 = \begin{bmatrix} I \\ -I \end{bmatrix}, \quad d_1 = \begin{bmatrix} \pi^{\text{nom}} + \rho \\ \rho - \pi^{\text{nom}} \end{bmatrix},$$

the DRKP problem becomes

$$\begin{array}{ll} \text{maximize} & \big(\min(\log(R^T b) + \lambda_+ - \lambda_-)\big) - \\ & (\pi^{\text{nom}})^T(\lambda_+ - \lambda_-) - \rho^T(\lambda_+ + \lambda_-) \\ \text{subject to} & b \in B, \quad \lambda_+ \geq 0, \quad \lambda_- \geq 0, \end{array}$$

with variables $b, \lambda_+, \lambda_-$. Defining $\lambda = \lambda_+ - \lambda_-$, we have $|\lambda| = \lambda_+ + \lambda_-$, so the DRKP becomes

$$\begin{array}{ll} \text{maximize} & \min(\log(R^T b) + \lambda) - (\pi^{\text{nom}})^T \lambda - \rho^T |\lambda| \\ \text{subject to} & b \in B, \end{array}$$

with variables $b, \lambda$. □

**Proof of Theorem 3**

*Proof of theorem* 3. We define $x = -\log(R^T b)$, $z = W^{-1}(\pi - \pi^{\text{nom}})$, and $D_{p,W} = \{z \mid \|z\|_p \leq 1,\ \mathbf{1}^T W z = 0,\ \pi^{\text{nom}} + Wz \geq 0\}$. Then we have

$$\begin{array}{rl} G_\Pi(b) & = -\sup_{\pi\in\Pi}((\pi - \pi^{\text{nom}})^T x + (\pi^{\text{nom}})^T x) \\ & = -\sup_{z\in D_{p,W}} z^T W^T x + (\pi^{\text{nom}})^T x \\ & = \sup_{\mu,\lambda\geq 0} \quad (-\sup_{\|z\|_p\leq 1} z^T W^T (x + \lambda - \mu 1) + \\ & \quad (\pi^{\text{nom}})^T(\lambda + x)) \\ & = \sup_{\mu,\lambda\geq 0} \quad (-\|W^T(x + \lambda - \mu 1)\|_q + \\ & \quad (\pi^{\text{nom}})^T(-\lambda - x)). \end{array}$$

Here the second last equation is the Lagrangian form where we keep the $p$-norm constraint as a convex indicator, and the last equation is based on the Hölder equality

$$\sup_{\|z\|_p\leq 1} z^T W^T(x + \lambda - \mu 1) = \|W^T(x + \lambda - \mu 1)\|_q,$$

Using Slater's condition, strong duality can easily be verified. Using this expression for $G_\Pi(b)$, and let $u = -x - \lambda = \log(R^T b) - \lambda \leq \log(R^T b)$, then the DCP formulation of DRKP becomes

$$\begin{array}{ll} \text{maximize} & (\pi^{\text{nom}})^T(u) - \|W^T(u - \mu 1)\|_q \\ \text{subject to} & u \leq \log(R^T b), \\ & b \in B, \end{array}$$

with variables $b, u, \mu$.

This DRKP problem follows DCP rule because $\pi^{\text{nom},T} u$ is a linear function of $u$, $-\|W^T(u - \mu 1)\|_q$ is a concave function of $u$ and $\mu$ for $q \geq 1$, and $\log(R^T b)$ is a concave function of $b$, hence $u \leq \log(R^T b)$ is a concave constraint. Therefore, this DRKP problem is maximizing a concave objective concave constraint problem that follows DCP rule. □

**Proof of Theorem** 4

*Proof of theorem* 4. We define $x = -\log(R^T b)$ again. Our goal is to minimize $-G_\Pi(b) = \sup_{\pi\in\Pi} \pi^T x$. We form a dual of this problem, working with the constraints $D_f(\pi\|\pi_0) \le \epsilon$ and $\mathbf{1}^T\pi = 1$; we keep the constraint $\pi \ge 0$ implicit. With dual variables $\lambda \in \mathbb{R}_+$, $\gamma \in \mathbb{R}$, then for $\pi \ge 0$, the Lagrangian is

$$\begin{aligned} L(\gamma,\lambda,\pi) &= \pi^T x + \lambda(-(\pi^{\text{nom}})^T f(\tfrac{\pi}{\pi^{\text{nom}}}) + \epsilon) \\ &\quad -\gamma(e^T\pi - 1) + I_+(\pi), \end{aligned}$$

where $I_+$ is the indicator function of $\mathbb{R}_+^K$. The dual objective function is

$$\begin{aligned} &\sup_{\pi\ge 0} L(\gamma,\lambda,\pi) \\ &= \sup_{\pi\ge 0}(\textstyle\sum_{i=1}^K \pi_i^{\text{nom}}(\frac{\pi_i}{\pi_i^{\text{nom}}}x_i - \frac{\pi_i}{\pi_i^{\text{nom}}}\gamma - \lambda f(\frac{\pi_i}{\pi_i^{\text{nom}}}))) \\ &\quad +\lambda\epsilon + \gamma \\ &= \textstyle\sum_{i=1}^K \pi_{0,i} \sup_{t_i\ge 0}(t_i(x_i-\gamma) - \lambda f(t_i)) + \lambda\epsilon + \gamma \\ &= \textstyle\sum_{i=1}^K \pi_{0,i}\lambda f^*(\frac{x_i-\gamma}{\lambda}) + \lambda\epsilon + \gamma. \end{aligned}$$

Using Slater's condition, strong duality can easily be verified. We can write the problem as

$$\begin{array}{ll} \text{maximize} & -(\pi^{\text{nom}})^T \lambda f^*(\frac{-\log(R^T b)-\gamma}{\lambda}) - \lambda\epsilon - \gamma \\ \text{subject to} & \lambda \ge 0, \quad b \in B, \end{array}$$

with variables $b, \gamma, \lambda$. We transform the problem to follow the disciplined convex programming (DCP) rules by convex relaxation of the equality constraint. Now DRKP becomes

$$\begin{array}{ll} \text{maximize} & -(\pi^{\text{nom}})^T w - \epsilon\lambda - \gamma \\ \text{subject to} & w \ge \lambda f^*(\frac{z}{\lambda}) \\ & z \ge -\log(R^T b) - \gamma \\ & \lambda \ge 0, \quad b \in B, \end{array}$$

with variables $b, \gamma, \lambda, w, z$.

Here

$$\lambda f^*(\frac{z}{\lambda}) = (\lambda f)^*(z) = \sup_{t\ge 0}(tz - \lambda f(t)),$$

is the perspective function of the non-decreasing convex function $f^*(z)$, so it is also a convex function that is non-decreasing in $z$. Additionally, $-\log(R^T b) - \gamma$ is a convex function of $b$ and $\gamma$; then from the DCP composition rule, we know this form of DRKP is convex. □

**Proof of Theorem** 5

*Proof of theorem* 5. The worst-case log growth $G_\Pi(b)$ is given by the value of the following LP,

$$\begin{array}{ll} \text{minimize} & \pi^T \log(R^T b) \\ \text{subject to} & Q\mathbf{1} = \pi, \quad Q^T\mathbf{1} = \pi^{\text{nom}}, \quad Q \ge 0, \\ & \sum_{i,j} Q_{ij}c_{ij} \le s, \end{array}$$

with variable $Q$. Using strong duality for LP, the DRKP becomes

$$\begin{array}{ll} \text{maximize} & \left(\sum_j \pi_j^{\text{nom}} \min_i(\log(R^T b)_i + \lambda c_{ij}) - s\lambda\right) \\ \text{subject to} & b \in B, \quad \lambda \ge 0. \end{array}$$

where $\lambda \in \mathbb{R}_+$ is the dual variable.

The problem follows the disciplined convex programming (DCP) rules, because $\log(R^T b)_i + \lambda c_{ij}$ is a concave function of $b$ and $\lambda$, therefore $\min_i(\log(R^T b)_i + \lambda c_{ij})$ is a concave function of $b$ and $\lambda$, then the entire objective is a concave function of $b$ and $\lambda$; and the constraint $b \in B$ and $\lambda \ge 0$ also follows DCP rule. □

**Proof of Theorem** 6

*Proof of theorem* 6. Following [28] lemma 1, to prove our theorem, we only need to verify Slater's constraint qualification conditions for strong duality and verify that $\log(r^T b)$ is integrable for all $\pi$ in $\Pi$.

For Slater's constraint qualification conditions, we only need to find a strictly feasible $\pi$ in $\Pi$. Since there exists $\pi_0 \in \Pi$ such that $\mathbf{E}_{\pi_0} r = \mu_0$, $\mathbf{E}_{\pi_0}(r-\mu_0)(r-\mu_0)^T = \sigma_0$, $\pi_0$ is a strictly feasible point in $\Pi$.

For integrablility, notice that $\mathbf{E}_\pi \log(r^T b)$ is on finite event space, so the result would naturally be finite for any $\Pi$ as a subset of $K$-dimensional probability simplex. □

## 8 Supplementary Material: Proof of theorems in section 3

**Proof of Theorem** 7 [Proof of theorem 7]

*Proof.* We have

$$\inf_{\pi_N \in \Pi} \mathbf{E}_{\pi_N}[\frac{S_N}{S_N^\star} \mid O_{N-1}] = \inf_{\pi_N \in \Pi} \mathbf{E}_{\pi_N}[\frac{V_N}{V_N^\star} \mid O_{N-1}]\frac{S_{N-1}}{S_{N-1}^\star}$$

If we can prove that for any $b_N$,

$$\inf_{\pi_N \in \Pi} \mathbf{E}_{\pi_N}(\frac{V_N}{V_N^\star} \mid O_{N-1}) \leq 1$$

then $\frac{S_N}{S_N^\star}$ is a decreasing semimartingale under non-linear expectation theory [49], with

$$\inf_{\pi_1,\ldots,\pi_{N-1},\pi_N \in \Pi} \mathbf{E}_{\pi_1,\ldots,\pi_{N-1},\pi_N} \lim_N \frac{S_N}{S_N^\star} \leq \inf_{\pi_1,\ldots,\pi_{N-1},\pi_{N-1} \in \Pi} \mathbf{E}_{\pi_1,\ldots,\pi_{N-1}} \lim_N \frac{S_{N-1}}{S_{N-1}^\star} \leq \frac{S_0}{S_0^\star} = 1$$

Now, for any $\epsilon > 0$, by the maximizing property of the distributional robust bet $b_N^\star$, we have

$$\inf_{\pi_N \in \Pi} \mathbf{E}(\epsilon \log(V_N) + (1-\epsilon)\log(V_N^\star) \mid O_{N-1}) \leq \inf_{\pi_N \in \Pi} \mathbf{E}(\log(V_N^\star) \mid O_{N-1})$$

Rewrite the left side, we have

$$\inf_{\pi_N \in \Pi} \mathbf{E}(\epsilon \log(V_N)+(1-\epsilon)\log(V_N^\star) \mid O_{N-1}) = \log(1-\epsilon) + \inf_{\pi_N \in \Pi} \mathbf{E}[\log(V_N^\star)+\log(1+\frac{\epsilon}{1-\epsilon}\frac{V_N}{V_N^\star}) \mid O_{N-1}]$$

Using superlinear property of $\mathcal{E} := \inf_{\pi_N \in \Pi} \mathbf{E}_{\pi_N}$, we have

$$\inf_{\pi_N \in \Pi} \mathbf{E}[\log(V_N^\star)+\log(1+\frac{\epsilon}{1-\epsilon}\frac{V_N}{V_N^\star}) \mid O_{N-1}] \geq \inf_{\pi_N \in \Pi} \mathbf{E}[\log(V_N^\star) \mid O_{N-1}] + \inf_{\pi_N \in \Pi} \mathbf{E}[\log(1+\frac{\epsilon}{1-\epsilon}\frac{V_N}{V_N^\star}) \mid O_{N-1}]$$

Therefore, combine this inequity with the previous inequity from the maximizing property of the distributional robust bet $b_N^\star$, we have

$$\inf_{\pi_N \in \Pi} \mathbf{E}[\frac{1}{\epsilon}\log(1+\frac{\epsilon}{1-\epsilon}\frac{V_N}{V_N^\star}) \mid O_{N-1}] \leq \frac{1}{\epsilon}\log(\frac{1}{1-\epsilon})$$

taking lower limit $\epsilon \to 0^+$, since $\frac{V_N}{V_N^\star}$) are bounded, from Fatou's lemma, we have

$$\begin{aligned}\inf_{\pi_N \in \Pi} \mathbf{E}[\frac{V_N}{V_N^\star}) \mid O_{N-1}] &= \inf_{\pi_N \in \Pi} \mathbf{E}[\underline{\lim}_{\epsilon \to 0^+} \frac{1}{\epsilon}\log(1+\frac{\epsilon}{1-\epsilon}\frac{V_N}{V_N^\star}) \mid O_{N-1}] \\ &\leq \underline{\lim}_{\epsilon \to 0^+} \inf_{\pi_N \in \Pi} \mathbf{E}[\frac{1}{\epsilon}\log(1+\frac{\epsilon}{1-\epsilon}\frac{V_N}{V_N^\star}) \mid O_{N-1}] \\ &\leq \underline{\lim}_{\epsilon \to 0^+} \frac{1}{\epsilon}\log(\frac{1}{1-\epsilon}) = 1\end{aligned}$$

□

**Proof of Theorem** 8

*Proof of theorem* 8. The two sequences

$$\frac{S_N}{S_N^\star} - \inf_{\pi_N\in\Pi} \mathbf{E}_{\pi_N}[\frac{S_N}{S_N^\star} \mid O_{N-1}] = \sum_{s=1}^{N}\{\log(V_s) - \log(V_s^\star) - \inf_{\pi_s\in\Pi} \mathbf{E}[\log(V_s) - \log(V_s^\star) \mid O_{s-1}]\}$$

and

$$\frac{S_N^\star}{S_N} - \inf_{\pi_N\in\Pi} \mathbf{E}_{\pi_s}[\frac{S_N^\star}{S_N} \mid O_{N-1}] = \sum_{s=1}^{N}\{\log(V_s^\star) - \log(V_s) - \inf_{\pi_s\in\Pi} \mathbf{E}[\log(V_s^\star) - \log(V_s) \mid O_{s-1}]\}$$

both form $G$-martingale sequences under nonlinear expectation $\mathcal{E} := \inf_{\pi_N\in\Pi} \mathbf{E}_{\pi_N}$ as defined in [49]. Therefore, from nonlinear expectation version of Doob's martingale convergence theorem [49, 50], both sequences converge to a finite value almost surely if $\sum_{s=1}^N\{\log(V_s^\star) - \log(V_s) - \inf_{\pi_s\in\Pi} \mathbf{E}[\log(V_s^\star) - \log(V_s) \mid O_{s-1}]\}$ is uniformly bounded almost surely for all $N$. From previous theorem, we know $\frac{S_N}{S_N^\star} \leq 1$. Therefore, $\overline{\lim}_{N\to\infty} \frac{S_N}{S_N^\star} = 0$ almost surely if and only if $\sum_{N=1}^\infty[\log(V_N) - \log(V_N^\star) - \inf_{\pi_N\in\Pi} \mathbf{E}_{\pi_N}[\log(V_N^\star) - \log(V_N) \mid O_{N-1}]] = \infty$ □

## 9 Details of the horse racing numerical example

In this section we provide additional details of the horse racing numerical example. Our example is a simple horse race with $n$ horses, with bets placed on each horse placing, *i.e.*, coming in first or second. There are thus $K = n(n-1)/2$ outcomes (indexed as $j, k$ with $j < k \leq n$), and $n$ bets (one for each horse to place). We first describe the nominal distribution of outcomes $\pi^{\text{nom}}$. We model the speed of the horses as independent random variables, with the fastest and second fastest horses placing. With this model, $\pi^{\text{nom}}$ is entirely described by the probability that horse $i$ comes in first, we which denote $\beta_i$. For $j < k$, we have

$$\begin{aligned}
\pi_{jk}^{\text{nom}} &= P(\text{horse } j \text{ and } k \text{ are the first two})\\
&= P(j \text{ is 1st}, k \text{ is 2nd}) + P(k \text{ is 1st}, j \text{ is 2nd})\\
&= P(j \text{ is 1st})P(k \text{ is 2nd} \mid j \text{ is 1st})\\
&\quad +P(k \text{ is 1st})P(j \text{ is 2nd} \mid k \text{ is 1st})\\
&= \beta_j(\beta_k/(1-\beta_j)) + \beta_k(\beta_j/(1-\beta_k))\\
&= \beta_j\beta_k(\tfrac{1}{1-\beta_i} + \tfrac{1}{1-\beta_j}).
\end{aligned}$$

The fourth line uses $P(k \text{ is 2nd} \mid j \text{ is 1st}) = \beta_k/(1-\beta_j)$.

For the return matrix, we use parimutuel betting, with the fraction of bets on each horse equal to $\beta_i$, the probability that it will win (under the nominal probability distribution). The return matrix $R \in \mathbb{R}^{n\times K}$ then has the form

$$R_{i,jk} = \begin{cases} \frac{n}{1+\beta_j/\beta_k} & i = j \\ \frac{n}{1+\beta_k/\beta_j} & i = k \\ 0 & i \notin \{j,k\}, \end{cases}$$

where we index the columns (outcomes) by the pair $jk$, with $j < k$.

Our set of possible distributions is the box

$$\Pi_\eta = \{\pi \mid |\pi - \pi^{\text{nom}}| \leq \eta\pi^{\text{nom}},\ \mathbf{1}^T\pi = 1,\ \pi \geq 0\},$$

where $\eta \in (0,1)$, *i.e.*, each probability can vary by $\eta$ from its nominal value.

Another uncertainty set is the ball

$$\Pi_c = \{\pi \mid \|\pi - \pi^{\text{nom}}\|_2 \leq c,\ \mathbf{1}^T\pi = 1,\ \pi \geq 0\}$$

For our specific example instance, we take $n = 20$ horses, so there are $K = 190$ outcomes. We choose $\beta_i$, the probability distribution of the winning horse, proportional to $\exp z_i$, where we sample independently $z_i \sim \mathcal{N}(0, 1/4)$. This results in $\beta_i$ ranging from around 20% (the fastest horse) to around 1% (the slowest horse).

First, we show growth rate and worst-case growth rate for the Kelly optimal and the distributional robust Kelly optimal bets under two uncertainty sets. In table 1 of the main paper, we show the comparison for box uncertainty set with $\eta = 0.26$ and for ball uncertainty set with $c = 0.016$. The two parameters are chosen so that the worst case growth of Kelly bets for both uncertainty sets are $-2.2\%$. In particular, using standard Kelly betting, we lose money (when the distribution is chosen as the worst one for the Kelly bets). We can see that, as expected, the Kelly optimal bet has higher log growth under the nominal distribution, and the distributional robust Kelly bet has better worst-case log growth. We see that the worst-case growth of the distributional robust Kelly bet is significantly better than the worst-case growth of the nominal Kelly optimal bet. In particular, with robust Kelly betting, we make money, even when the worst distribution is chosen. The nominal Kelly optimal bet $b^{\mathrm{K}}$ and the distributional robust Kelly bet $b^{\mathrm{RK}}$ for both uncertainty sets in figure 1 of main paper. For each of our bets $b^{\mathrm{K}}$ and $b^{\mathrm{RK}}$ shown above, we find a corresponding worst case distribution, denoted $\pi^{\mathrm{wc,K}}$ and $\pi^{\mathrm{wc,RK}}$, which minimize $G_\pi(b)$ over $\pi \in \Pi$. These distributions, shown for box uncertainty set and ball uncertainty set in figure 3 of main paper, achieve the corresponding worst-case log growth for the two bet vectors. Finally, we compare the expected wealth logarithmic growth rate as we increase the size of the uncertainty sets. For the box uncertainty set we choose $\eta \in [0, 0.3]$, and for the ball uncertainty set we choose $c \in [0, 0.02]$, we look at the expected growth for both the kelly bet $b^K$ and the distributional robust kelly bet $b^{RK}$ under both the nominal probability $\pi^{\mathrm{nom}}$ and the worst case probability $\pi^{\mathrm{worst}}$.

# 10 Supplementary Material: Details for learning uncertainty set through bi-level optimization

To use distributional robust Kelly strategy, one of the limitation is the requirement to tune of the parameters for the uncertainty set. Tuning is often done by hand, or by simple methods such as a crude grid search. In this section we propose a method to automate this process, by adjusting the parameters using an approximate gradient of the performance metric with respect to the parameters.

To give more colors to the general setup presented in the discussion section (section 5) of main paper, we consider a more concrete setting of learning uncertainty with features (covariates). Assuming that we are playing sequential gambling games, where each game is characterized by features (covariates) $X$. We represent the nominal return distribution as a function of features: $\pi_0 = h_{\theta_0}(X)$, here $h_{\theta_0}$ is a logistic function $h_{\theta_0}(X) = \exp(\theta_0^T X)/Z$ or in general a deep neural network with parameter $\theta_0$, the radius/shape of the uncertainty set is parametrized by $\theta_1$. For example, in the transformed $l_p$ ball uncertainty set, $\theta_1 = W$. The uncertainty set $\Pi(\theta; X)$ has parameter $\theta = (\theta_0, \theta_1)$.

For a given uncertainty set with parameter $\theta$ and feature $X$, the distributional robust Kelly strategy is the solution to the convex optimization problem:

$$b^\star(\theta; X) = \arg\max_{b \in B} \quad \min_{\pi \in \Pi(\theta; X)} \mathbf{E}_\pi \log(r^T b).$$

To automatically tuning the parameters, we could define a performance metric with respect to the parameters $L(\theta)$. For example, given out of samples observations $X_i$, $r_i$ that are not accessible when we fit the mean the functional form of $\Pi$, we could define an out-of-sample Kelly loss as the performance metric:

$$\max_\theta \quad L(\theta) = \frac{1}{N} \sum_{i=1}^{N} \log(r_i^T b^\star(\theta; X_i))$$

Using this performance metric as the high level objective, we could choose the uncertainty set parameter $\theta$ through bi-level optimization using approximate gradient method. Our method relies on recently developed methods[58] that can efficiently evaluate the derivative of the solution of a disciplined convex optimization problem with respect to its parameters. Using previous result in this paper, this distributional robust optimization problem can be transformed into a disciplined convex optimization, which allows automated differentiation with respect to the solution map $b^\star(\theta; X)$. We have provided CVXPY code for the uncertainty sets in this paper, which would allow users to use them through the recently developed software framework to embed our problem into differential programming framework like Tensorflow and PyTorch to learn the uncertainty set via deep learning.

## 11 Supplementary Material: CVXPY example codes

All of the formulations of distributional robust Kelly problem (DRKP) are not only tractable, but easily expressed in domain specific language for convex optimization. The CVXPY code to specify and solve the DRKP for ball and box constraints, for example, is given below.

For box uncertainty set, $\Pi_\rho = \{\pi \mid |\pi - \pi^{\text{nom}}| \le \rho,\ \mathbf{1}^T\pi = 1,\ \pi \ge 0\}$, the CVXPY code is

```
pi_nom = Parameter(K, nonneg=True)
rho = Parameter(K,nonneg=True)
b = Variable(n)
mu = Variable(K)
wc_growth_rate = min(log(R.T*b) + mu)
                 -pi_nom.T*abs(mu)
                 -rho.T*mu
constraints = [sum(b) == 1, b >= 0]
DRKP = Problem(Maximize(wc_growth_rate),
               constraints)
DRKP.solve()
```

For ball uncertainty set, $\Pi_c = \{\pi \mid \|\pi - \pi^{\text{nom}}\|_2 \le c,\ \mathbf{1}^T\pi = 1,\ \pi \ge 0\}$, the CVXPY code is

```
pi_nom = Parameter(K, nonneg=True)
c = Parameter((1,1), nonneg=True)
b = Variable(n)
U = Variable(K)
mu = Variable(K)
log_growth = log(R.T*b)
wc_growth_rate = pi_nom.T*F-c*norm(U-mu,2)
constraints = [sum(b) == 1,
               b >= 0,
               U <= log_growth]
DRKP = Problem(Maximize(wc_growth_rate),
               constraints)
DRKP.solve()
```

Here `R` is the matrix whose columns are the return vectors, `pi_nom` is the vector of nominal probabilities. `rho` is $K$ dimensional box constraint and `c` is radius of the ball. For each problem, the second to last line forms the problem, and in the last line the problem is solved. The robust optimal bet is written into `b.value`.

To learn uncertainty set, we could parametrize $\pi_0$ via logistic model $\pi_0 = \text{softmax}(\theta X)$, where $\theta \in \mathbb{R}^{K\times M}$ is the parameter to learn, $X \in \mathbb{R}^M$ is a feature of each game.

The full python notebook code for the horse gambling example is also attached at appendix. The computational resource used is fairly light-weight. All the computation in the notebook is done via Google's Colab, and the notebook is also easy to run with laptops.

To give a taste of the framework, for ball uncertainty set, $\Pi_\rho = \{\pi \mid |\pi - \pi^{\text{nom}}| \le \rho,\ \mathbf{1}^T\pi = 1,\ \pi \ge 0\}$, the CVXPY code to build CvxpyLayer is

```
import cvxpy as cvx
import torch
from cvxpylayers.torch import CvxpyLayer
# generate_ball_problem
pi_0 = cvx.Parameter(K, nonneg=True)
rho = cvx.Parameter(K, nonneg=True)
R_cvx = cvx.Parameter((n,K), nonneg=True)
b = cvx.Variable(n)
mu = cvx.Variable(K)
log_growth = cvx.log(R_cvx.T@b )
rob_growth_rate = cvx.min(log_growth + mu  )
rob_growth_rate = rob_growth_rate - pi_0.T@mu - rho.T@cvx.abs(mu)
constraints = [cvx.sum(b) == 1, b >= 0]
DRKP = cvx.Problem(cvx.Maximize(rob_growth_rate), constraints)
```

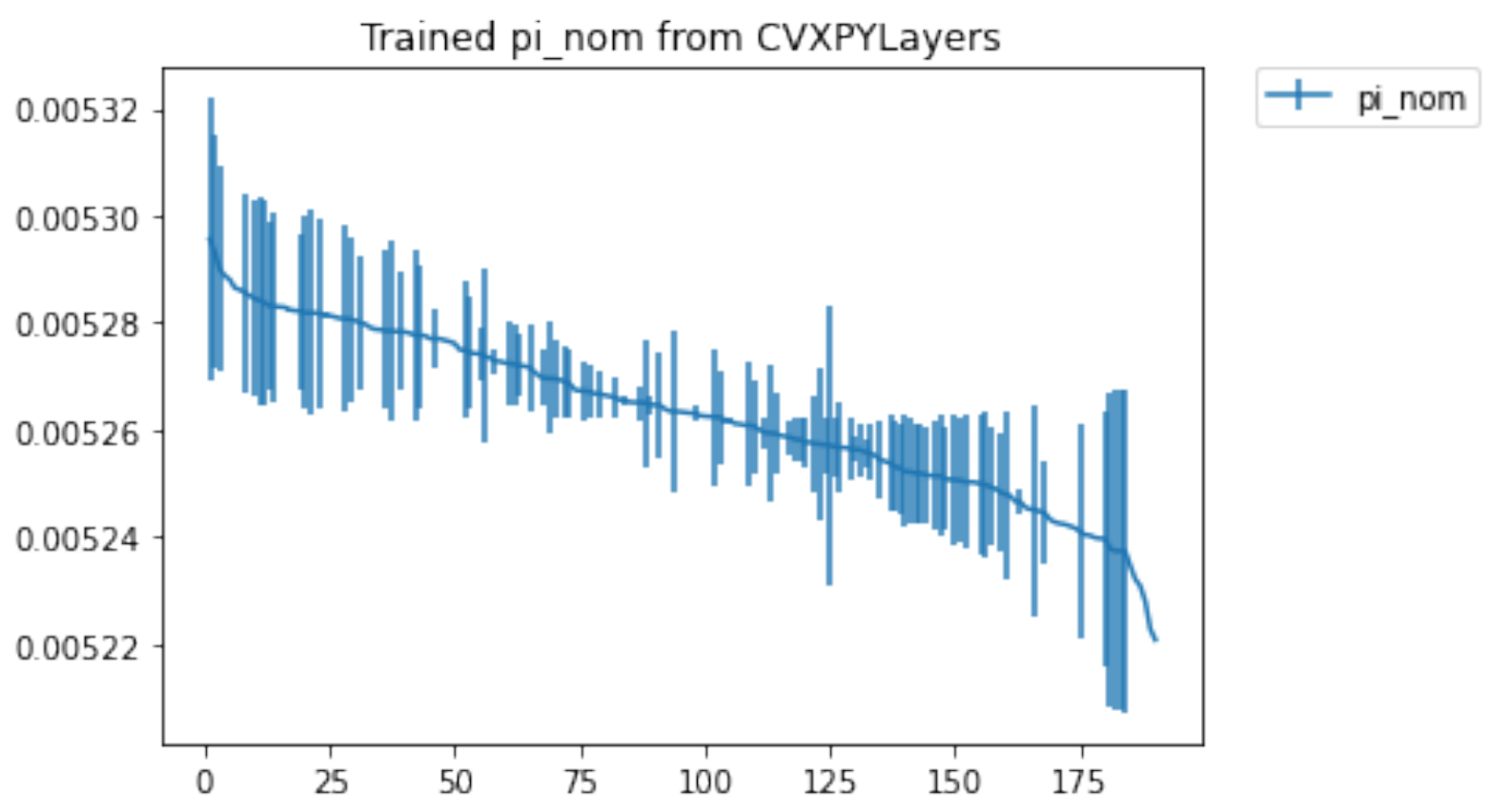

**Figure 4:** The trained $\pi_{\text{nom}}$ with error bar from box constraint, $\pi_{\text{nom}}$ is initialized at uniform distribution.

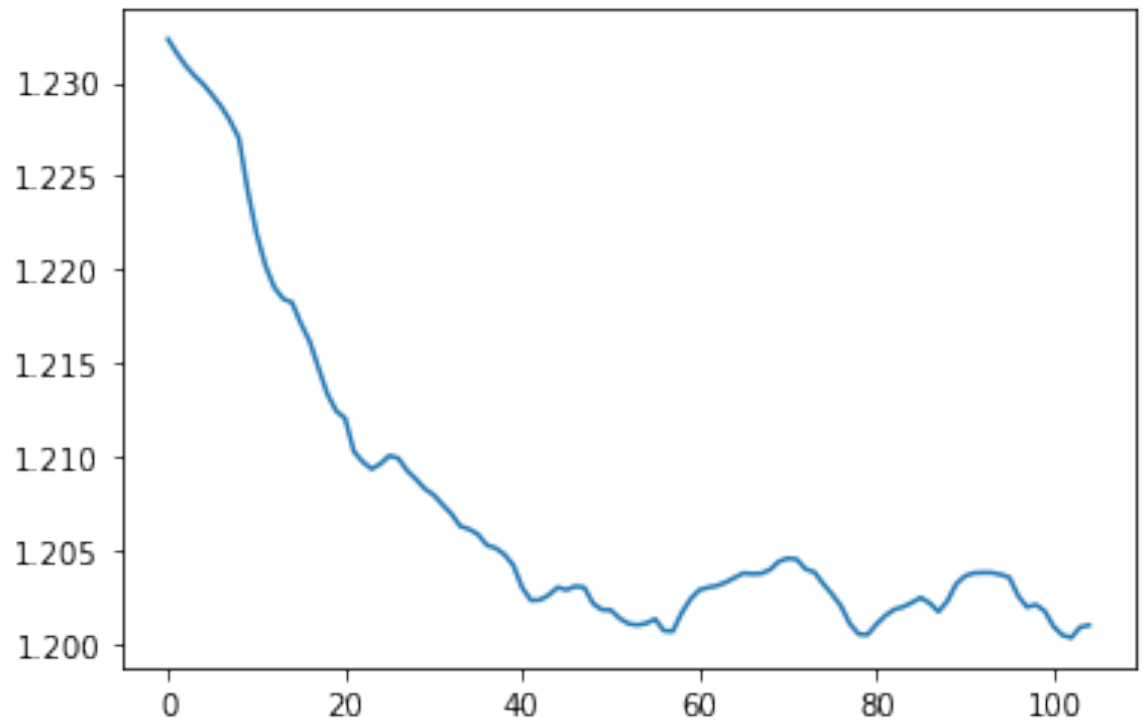

**Figure 5:** The training loss for training $\rho$ using projected ADAM optimizer (onto non-negative vectors), $\rho$ is initialized at $rho = (0.1, \dots, 0.1)$ using PyTorch with initial step size $10^{-6}$.

```
problem = DRKP
parameters=[R_cvx, pi_0, rho]
policy =  CvxpyLayer(problem, parameters, [b])
```

The training code using PyTorch looks like:

```
# Initialize:
Rho = torch.from_numpy(np.ones(K)*1e-4).requires_grad_(True)
log_Pi_0 = torch.from_numpy(np.zeros(K)).requires_grad_(True)
torch_variables = [log_Pi_0, Rho]
R_torch = torch.from_numpy(R)
Pi_test_torch = torch.from_numpy(Pi_test)
# Loss:
def evaluate( R_torch, log_Pi_0, Rho, Pi_test_torch):
    Pi_0 = torch.nn.functional.softmax(log_Pi_0)
    b, = policy(R_torch, Pi_0, Rho)
    logs = torch.log(R_torch.T @ b)
    cost = - torch.sum(Pi_test_torch*logs[None,:])
    return   cost

# Training:
iters = 100
results = []
optimizer = torch.optim.Adam(torch_variables, lr=1e-2)
for i in range(iters):
    optimizer.zero_grad()
```

```
loss = evaluate(R_torch, log_Pi_0, Rho, Pi_test_torch)
loss.backward()
optimizer.step()
# Project so that Rho is non-negative
Rho.data = torch.max(Rho.data,torch.zeros_like(Rho.data))
results.append(loss.item())
print("(iter %d) loss: %g " % (i, results[-1]))
```